%Paper: alg-geom/9210004
%From: debarre@msri.org (Olivier Debarre)
%Date: Tue, 13 Oct 92 12:27:49 PDT

\magnification=1200

\hsize=13.2cm
\vsize=19.5cm
\parindent=0cm   \parskip=0pt
\pageno=1

\def\ind{\hskip 1cm\relax}
\hoffset=15mm
\voffset=1cm

\ifnum\mag=\magstep1
\hoffset=0.2cm   % offset horizontal en \magnification=1200
\voffset=-.5cm   % offset vertical en \magnification=1200
\fi

\pretolerance=500 \tolerance=1000  \brokenpenalty=5000

\catcode`\@=11

\font\eightrm=cmr8         \font\eighti=cmmi8
\font\eightsy=cmsy8        \font\eightbf=cmbx8
\font\eighttt=cmtt8        \font\eightit=cmti8
\font\eightsl=cmsl8        \font\sixrm=cmr6
\font\sixi=cmmi6           \font\sixsy=cmsy6
\font\sixbf=cmbx6

\skewchar\eighti='177 \skewchar\sixi='177
\skewchar\eightsy='60 \skewchar\sixsy='60

\newfam\gothfam           \newfam\bboardfam
\newfam\cyrfam

\def\tenpoint{%
  \textfont0=\tenrm \scriptfont0=\sevenrm \scriptscriptfont0=\fiverm
  \def\rm{\fam\z@\tenrm}%
  \textfont1=\teni  \scriptfont1=\seveni  \scriptscriptfont1=\fivei
  \def\oldstyle{\fam\@ne\teni}\let\old=\oldstyle
  \textfont2=\tensy \scriptfont2=\sevensy \scriptscriptfont2=\fivesy
  \textfont\gothfam=\tengoth \scriptfont\gothfam=\sevengoth
  \scriptscriptfont\gothfam=\fivegoth
  \def\goth{\fam\gothfam\tengoth}%
  \textfont\bboardfam=\tenbboard \scriptfont\bboardfam=\sevenbboard
  \scriptscriptfont\bboardfam=\sevenbboard
  \def\bb{\fam\bboardfam\tenbboard}%
 \textfont\cyrfam=\tencyr \scriptfont\cyrfam=\sevencyr
  \scriptscriptfont\cyrfam=\sixcyr
  \def\cyr{\fam\cyrfam\tencyr}%
  \textfont\itfam=\tenit
  \def\it{\fam\itfam\tenit}%
  \textfont\slfam=\tensl
  \def\sl{\fam\slfam\tensl}%
  \textfont\bffam=\tenbf \scriptfont\bffam=\sevenbf
  \scriptscriptfont\bffam=\fivebf
  \def\bf{\fam\bffam\tenbf}%
  \textfont\ttfam=\tentt
  \def\tt{\fam\ttfam\tentt}%
  \abovedisplayskip=12pt plus 3pt minus 9pt
  \belowdisplayskip=\abovedisplayskip
  \abovedisplayshortskip=0pt plus 3pt
  \belowdisplayshortskip=4pt plus 3pt
  \smallskipamount=3pt plus 1pt minus 1pt
  \medskipamount=6pt plus 2pt minus 2pt
  \bigskipamount=12pt plus 4pt minus 4pt
  \normalbaselineskip=12pt
  \setbox\strutbox=\hbox{\vrule height8.5pt depth3.5pt width0pt}%
  \let\bigf@nt=\tenrm       \let\smallf@nt=\sevenrm
  \normalbaselines\rm}

\def\eightpoint{%
  \textfont0=\eightrm \scriptfont0=\sixrm \scriptscriptfont0=\fiverm
  \def\rm{\fam\z@\eightrm}%
  \textfont1=\eighti  \scriptfont1=\sixi  \scriptscriptfont1=\fivei
  \def\oldstyle{\fam\@ne\eighti}\let\old=\oldstyle
  \textfont2=\eightsy \scriptfont2=\sixsy \scriptscriptfont2=\fivesy
  \textfont\gothfam=\eightgoth \scriptfont\gothfam=\sixgoth
  \scriptscriptfont\gothfam=\fivegoth
  \def\goth{\fam\gothfam\eightgoth}%
  \textfont\cyrfam=\eightcyr \scriptfont\cyrfam=\sixcyr
  \scriptscriptfont\cyrfam=\sixcyr
  \def\cyr{\fam\cyrfam\eightcyr}%
  \textfont\bboardfam=\eightbboard \scriptfont\bboardfam=\sevenbboard
  \scriptscriptfont\bboardfam=\sevenbboard
  \def\bb{\fam\bboardfam}%
  \textfont\itfam=\eightit
  \def\it{\fam\itfam\eightit}%
  \textfont\slfam=\eightsl
  \def\sl{\fam\slfam\eightsl}%
  \textfont\bffam=\eightbf \scriptfont\bffam=\sixbf
  \scriptscriptfont\bffam=\fivebf
  \def\bf{\fam\bffam\eightbf}%
  \textfont\ttfam=\eighttt
  \def\tt{\fam\ttfam\eighttt}%
  \abovedisplayskip=9pt plus 3pt minus 9pt
  \belowdisplayskip=\abovedisplayskip
  \abovedisplayshortskip=0pt plus 3pt
  \belowdisplayshortskip=3pt plus 3pt
  \smallskipamount=2pt plus 1pt minus 1pt
  \medskipamount=4pt plus 2pt minus 1pt
  \bigskipamount=9pt plus 3pt minus 3pt
  \normalbaselineskip=9pt
  \setbox\strutbox=\hbox{\vrule height7pt depth2pt width0pt}%
  \let\bigf@nt=\eightrm     \let\smallf@nt=\sixrm
  \normalbaselines\rm}

\def\pc#1{\bigf@nt#1\smallf@nt}         \def\pd#1 {{\pc#1} }

\def\^#1{\if#1i{\accent"5E\i}\else{\accent"5E #1}\fi}
\def\"#1{\if#1i{\accent"7F\i}\else{\accent"7F #1}\fi}

\newtoks\auteurcourant      \auteurcourant={\hfil}
\newtoks\titrecourant       \titrecourant={\hfil}

\newtoks\hautpagetitre      \hautpagetitre={\hfil}
\newtoks\baspagetitre       \baspagetitre={\hfil}

\newtoks\hautpagegauche
\hautpagegauche={\eightpoint\rlap{\folio}\hfil\the\auteurcourant\hfil}
\newtoks\hautpagedroite
\hautpagedroite={\eightpoint\hfil\the\titrecourant\hfil\llap{\folio}}

\newtoks\baspagegauche      \baspagegauche={\hfil}
\newtoks\baspagedroite      \baspagedroite={\hfil}

\newif\ifpagetitre          \pagetitretrue

\footline={\ifpagetitre\the\baspagetitre\else
\ifodd\pageno\the\baspagedroite\else\the\baspagegauche\fi\fi
\global\pagetitrefalse}

\def\raggedbottom{\topskip 10pt plus 36pt\r@ggedbottomtrue}

\def\pointir{\unskip . --- \ignorespaces}

\def\Bigbreak{\vskip-\lastskip\bigbreak}
\def\Medbreak{\vskip-\lastskip\medbreak}

\def\ctexte#1\endctexte{%
  \hbox{$\vcenter{\halign{\hfill##\hfill\crcr#1\crcr}}$}}

\long\def\ctitre#1\endctitre{%
    \ifdim\lastskip<24pt\vskip-\lastskip\bigbreak\bigbreak\fi
  		\vbox{\parindent=0pt\leftskip=0pt plus 1fill
          \rightskip=\leftskip
          \parfillskip=0pt\bf#1\par}
    \bigskip\nobreak}

\long\def\section#1\endsection{%
\vskip 0pt plus 3\normalbaselineskip
\penalty-250
\vskip 0pt plus -3\normalbaselineskip
\Bigbreak
\message{[section \string: #1]}{\bf#1\unskip}\pointir}

\long\def\sectiona#1\endsection{%
\vskip 0pt plus 3\normalbaselineskip
\penalty-250
\vskip 0pt plus -3\normalbaselineskip
\Bigbreak
\message{[sectiona \string: #1]}%
{\bf#1}\medskip\nobreak}

\long\def\subsection#1\endsubsection{%
\Medbreak
{\it#1\unskip}\pointir}

\long\def\subsectiona#1\endsubsection{%
\Medbreak
{\it#1}\par\nobreak}

\def\rem#1\endrem{%
\Medbreak
{\it#1\unskip} : }

\def\remp#1\endrem{%
\Medbreak
{\pc #1\unskip}\pointir}

\def\rema#1\endrem{%
\Medbreak
{\it #1}\par\nobreak}

\def\newparwithcolon#1\endnewparwithcolon{
\Medbreak
{#1\unskip} : }

\def\newparwithpointir#1\endnewparwithpointir{
\Medbreak
{#1\unskip}\pointir}

\def\newpara#1\endnewpar{
\Medbreak
{#1\unskip}\smallskip\nobreak}

\let\+=\tabalign

\def\signature#1\endsignature{\vskip 15mm minus 5mm\rightline{\vtop{#1}}}

\mathcode`A="7041 \mathcode`B="7042 \mathcode`C="7043 \mathcode`D="7044
\mathcode`E="7045 \mathcode`F="7046 \mathcode`G="7047 \mathcode`H="7048
\mathcode`I="7049 \mathcode`J="704A \mathcode`K="704B \mathcode`L="704C
\mathcode`M="704D \mathcode`N="704E \mathcode`O="704F \mathcode`P="7050
\mathcode`Q="7051 \mathcode`R="7052 \mathcode`S="7053 \mathcode`T="7054
\mathcode`U="7055 \mathcode`V="7056 \mathcode`W="7057 \mathcode`X="7058
\mathcode`Y="7059 \mathcode`Z="705A

\def\spacedmath#1{\def\packedmath##1${\bgroup\mathsurround=0pt ##1\egroup$}%
\mathsurround#1 \everymath={\packedmath}\everydisplay={\mathsurround=0pt }}

\def\nospacedmath{\mathsurround=0pt \everymath={}\everydisplay={} }

\long\def\th#1 #2\enonce#3\endth{%
   \Medbreak
   {\pc#1} {#2\unskip}\pointir{\it #3}\medskip}

\long\def\tha#1 #2\enonce#3\endth{%
   \Medbreak
   {\pc#1} {#2\unskip}\par\nobreak{\it #3}\medskip}

\long\def\Th#1 #2 #3\enonce#4\endth{%
   \Medbreak
   #1 {\pc#2} {#3\unskip}\pointir{\it #4}\medskip}

\long\def\Tha#1 #2 #3\enonce#4\endth{%
   \Medbreak
   #1 {\pc#2} #3\par\nobreak{\it #4}\medskip}

\def\decale#1{\smallbreak\hskip 28pt\llap{#1}\kern 5pt}
\def\decaledecale#1{\smallbreak\hskip 34pt\llap{#1}\kern 5pt}
\def\puce{\smallbreak\hskip 6pt{$\scriptstyle\bullet$}\kern 5pt}

\def\displaylinesno#1{\displ@y\halign{
\hbox to\displaywidth{$\@lign\hfil\displaystyle##\hfil$}&
\llap{$##$}\crcr#1\crcr}}

\def\ldisplaylinesno#1{\displ@y\halign{
\hbox to\displaywidth{$\@lign\hfil\displaystyle##\hfil$}&
\kern-\displaywidth\rlap{$##$}\tabskip\displaywidth\crcr#1\crcr}}

\def\eqalign#1{\null\,\vcenter{\openup\jot\m@th\ialign{
\strut\hfil$\displaystyle{##}$&$\displaystyle{{}##}$\hfil
&&\quad\strut\hfil$\displaystyle{##}$&$\displaystyle{{}##}$\hfil
\crcr#1\crcr}}\,}

\def\system#1{\left\{\null\,\vcenter{\openup1\jot\m@th
\ialign{\strut$##$&\hfil$##$&$##$\hfil&&
        \enskip$##$\enskip&\hfil$##$&$##$\hfil\crcr#1\crcr}}\right.}

\let\@ldmessage=\message

\def\message#1{{\def\pc{\string\pc\space}%
                \def\'{\string'}\def\`{\string`}%
                \def\^{\string^}\def\"{\string"}%
                \@ldmessage{#1}}}

\def\up#1{\raise 1ex\hbox{\smallf@nt#1}}

\def\qed{\raise -2pt\hbox{\vrule\vbox to 10pt{\hrule width 4pt
                 \vfill\hrule}\vrule}}

\def\cqfd{\unskip\penalty 500\quad\vrule height 4pt depth 0pt width
4pt\medbreak}

\def\virg{\raise .4ex\hbox{,}}   % virgule aprs une fraction

 % point-virgule de ponctuation en maths

\def\build#1_#2^#3{\mathrel{
\mathop{\kern 0pt#1}\limits_{#2}^{#3}}}

\def\boxit#1#2{%
\setbox1=\hbox{\kern#1{#2}\kern#1}%
\dimen1=\ht1 \advance\dimen1 by #1 \dimen2=\dp1 \advance\dimen2 by #1
\setbox1=\hbox{\vrule height\dimen1 depth\dimen2\box1\vrule}%
\setbox1=\vbox{\hrule\box1\hrule}%
\advance\dimen1 by .4pt \ht1=\dimen1
\advance\dimen2 by .4pt \dp1=\dimen2  \box1\relax}

\def\date{\the\day\ \ifcase\month\or janvier\or f\'evrier\or mars\or
avril\or mai\or juin\or juillet\or ao\^ut\or septembre\or octobre\or
novembre\or d\'ecembre\fi \ {\old \the\year}}

\def\crog{{\vrule height 2.57mm depth 0.85mm width 0.3mm}\kern -0.36mm
[}

\def\crod{]\kern -0.4mm{\vrule height 2.57mm depth 0.85mm
width 0.3 mm}}

\def\rond{\kern 1pt{\scriptstyle\circ}\kern 1pt}

\def\hfl#1#2{\nospacedmath\smash{\mathop{\hbox to
12mm{\rightarrowfill}}\limits^{\scriptstyle#1}_{\scriptstyle#2}}}

\def\phfl#1#2{\nospacedmath\smash{\mathop{\hbox to
8mm{\rightarrowfill}}\limits^{\scriptstyle#1}_{\scriptstyle#2}}}

\catcode`\@=12

\showboxbreadth=-1  \showboxdepth=-1
\baselineskip=14pt
\spacedmath{1.7pt}
\baspagegauche={\centerline{\tenbf\folio}}
\baspagedroite={\centerline{\tenbf\folio}}
\hautpagegauche={\hfil}
\hautpagedroite={\hfil}
\font\eightrm=cmr10 at 8pt
\def\pa{\S\kern.15em}
\def\ra{\rightarrow}

\def\saut{\vskip 5mm plus 1mm minus 2mm}
\font\pc=cmcsc10 \rm
\def\Z{\hbox{\bf Z}}
\def\Q{\hbox{\bf Q}}

\def\Div{\mathop{\rm Div}\nolimits}

\def\Cl{\mathop{\rm Cl}\nolimits}

\def\rank{\mathop{\rm rank}\nolimits}

\parskip=2mm
\def\ie{\hbox{i.e. }}
\def\page{\hbox{p. }}

\def\bpf{base-point-free }

\saut
\ctitre
{\bf POINTS OF LOW DEGREE ON SMOOTH PLANE CURVES}
\endctitre
\medskip
\centerline{{\pc Olivier Debarre}
\footnote{(*)}{\rm Partially supported by the European Science
Project ``Geometry of Algebraic Varieties", Contract no. SCI-0398-C (A) and by
N.S.F. Grant DMS 92-03919.}
 and {\pc Matthew J. Klassen}}

\vskip 1cm

{\bf 1. Introduction}

\ind The purpose of this note is to provide some applications of a theorem of
Faltings ([Fa1]) to smooth plane curves, using ideas from [A] and [AH].

\ind Let $C$ be a smooth projective plane curve defined by an
equation of degree $d$ with rational coefficients. We show:

{\pc Theorem} 1.-- {\it
If $d\ge 7$, the curve $C$  has only
        finitely many points whose field of definition has degree
        $\leq d-2$ over $\Q$.}

\ind The result still holds for $d<7$, provided that the complex curve $C$ has
no morphisms of
degree $\le d-2$ onto an elliptic curve, an assumption which we show
automatically holds for $d\ge
7$. This result is sharp in the sense that if $C$ has a rational point, there
exist infinitely many
points on $C$ with field of definition of degree $\le d-1$. These points come
from the
intersection of $C$ with a rational line through a rational point. We show
further:

{\pc Theorem} 2.-- {\it If $d\ge 8$, all but finitely many points of  $C$
       whose field of definition has degree
        $ \le d-1$ over $\Q$ arise as the intersection of $C$ with a rational
line
through a rational point of $C$.}

\ind  In particular, if $C$ has no rational points, there are only finitely
many points whose field of definition has degree
        $\leq d-1$ over $\Q$.

\ind Again,
the result still holds for $d=6$ or $7$, provided that $C$ has no morphisms of
degree $\le d-1$
onto an elliptic curve, and for $d=5$, provided that $C$ has no morphisms  onto
an elliptic curve.

\ind Both results remain valid if $\Q $ is replaced by any number field.

\ind These theorems apply in particular to the Fermat curves $F_d$ with
equation
$X^d+Y^d=Z^d$, which is the case we had in mind when we started this
investigation.
Moreover, we can extend the results to all $d\ge 3$ in this case, with the one
exception
$d\ne 6$ (see \S6).

\ind The first author would like to thank the M.S.R.I., where  this work was
done,
for its support and hospitality.
\medskip
{\bf 2. Notation}

\ind For a  projective curve $C$, we denote by $C^{(d)}$ the symmetric product
varieties of $C$. We also denote by $J(C)$ the
Jacobian variety of $C$ and by   $W_d(C)$ the image of $C^{(d)}$ under the
Abel-Jacobi map to
$J(C)$ defined with respect to a chosen base point on $C$. It corresponds to
isomorphism
classes of line bundles on $C$ of degree $d$ which have a non-zero section.

\medskip
{\bf 3. Faltings' theorem}

\ind We first remark that
 theorem 1 is equivalent to the statement that the set  of
$\Q$--rational points  on the symmetric product $C^{(n)}$ is finite for any
$n\le d-2$.  This is
obtained by simply observing that any point on $C(K)$ with $[K:\Q ] = n$,
 together with its conjugates, forms a divisor of degree $n$ invariant by
${\rm Gal}\,({\overline{\Q}} / \Q )$, and hence defines a $\Q$--rational point
on
$C^{(n)}$. Furthermore, since $C$
is a smooth plane curve of degree $d$, it has no pencils of degree $\leq d-2$
([ACGH] \page 56,
exercise 18.(i)), hence $C^{(n)}$
maps isomorphically onto $W_n(C)$. Thus the proof of theorem 1 is reduced to
showing that
$W_n(C)(\Q )$ is finite for all $n\le d-2$. We use
the following beautiful result of Faltings:

{\pc Theorem} (Faltings, [Fa1]) -- {\it  Let $A$ be an abelian variety defined
over a number
field $K$.
        If $X$ is a subvariety of $A$ which does not contain any translate of a
positive-dimensional abelian subvariety of $A$, then $X$ contains only
finitely many $K$--rational points.}

\ind It is therefore enough to show that $W_n(C)$ does not
contain any non-zero abelian variety.

\ind The situation in theorem 2 is a bit more complicated since the morphism:
$$
\psi : C^{(d-1)}\longrightarrow W_{d-1}(C)$$
is no longer an isomorphism: each pencil of degree $d-1$ on $C$ corresponds to
a rational curve
in $C^{(d-1)}$ which is contracted by $\psi$. By [ACGH] \page 56, exercise
18.(ii), all such
pencils are given by the lines through a fixed point $x$ of $C$. Let $R_x$ be
the corresponding
rational curve in $C^{(d-1)}$. Since $\psi$ induces an isomorphism outside of
the union of all
$R_x$, any rational point of
$C^{(d-1)}$ corresponds either to a rational point of $W_{d-1}(C)$, or to a
rational point of
some $R_x$. Now let $x$ be a point of $C$ such that $R_x(\Q )$ is non-empty and
let $D$ be
the divisor on $C$ that corresponds to a point of $R_x(\Q )$. The points of $D$
are then on
a unique line $l$ (which passes through $x$) and, since $D$ is invariant under
the action of
${\rm Gal}\,({\overline{\Q}} / \Q )$, so is $l$, which is therefore rational.
It follows
that the divisor $l\cdot C$ is rational, hence so is $x=l\cdot C -D$. This
reduces the
proof of theorem 2  to showing that $W_{d-1}(C)(\Q )$ is finite. As above, by
Faltings'
theorem, it is enough to show that $W_{d-1}(C)$ does not contain any
positive-dimensional
abelian variety. \medskip
{\bf 4. Linear systems on smooth plane curves}

\ind Before proceeding to the proof of the  theorems, we gather here some
elementary  facts about
linear systems on smooth plane curves, which we will deduce from the following
result of Coppens
and Kato. Let $H$ be a hyperplane section on a smooth plane curve $C$ of degree
$d$ and let
$D$ be an effective divisor on $C$  which belongs to a \bpf pencil. Then we
have:

{\pc Theorem} (Coppens-Kato, [CK]) -- {\it  If $n<k(d-k)$ for some integer $k$,
the linear
system $|(k-1)H-D|$ is non-empty.}

\ind We assume now that $d\ge 5$. Here are the consequences that we need:

(4.1) {\it If $\deg (D)\le 2d-5$,  then either $D\equiv H$ or $D\equiv
H-x$ for some point $x$  on $C$.}

(4.2) {\it If $\deg (D)= 2d-4$,  then   $D\equiv 2H-x_1-x_2-x_3-x_4$ for some
points $x_1$, $x_2$, $x_3$ and $x_4$  on $C$, no three of them collinear.}

This follows from the theorem with $k=3$, except for $d=5$. In the latter case,
Riemmann-Roch says
that the $6$ points of $D$ are on a  conic, which is what we need.

(4.3) {\it If $\deg (D)= 2d-3$ and $d\ge 7$,  then $D\equiv 2H-x_1-x_2-x_3$ for
some points $x_1$,
$x_2$ and $x_3$  on $C$, not collinear.}

(4.4) {\it If $\deg (D)= 2d-2$ and $d\ge 8$,   then $D\equiv
2H-x_1-x_2$ for some points $x_1$ and $x_2$  on $C$.}

(4.5) {\it If $\deg (D)= 2d-2$, $\dim |D|\ge 2$ and $d\ge 6$,  then $D\equiv
2H-x_1-x_2$ for some
points $x_1$ and $x_2$  on $C$. In particular $\dim |D|=3$.}

For $d\ge 7$, the linear system $|2H-(D-x)|$ is non-empty for any point $x$ on
$C$ by (4.3). By
Riemann-Roch, this implies that $|2H-D|$ is non-empty. For $d=6$, by
Riemann-Roch, the 10 points
of $D$ are on $3$ linearly independent cubics, which must be reducible. Since
$7$
points of $D$ cannot be on a line (because $d=6$), all points of $D$ are on a
conic.

(4.6) {\it If $\deg (D)= 2d-2$ and $\dim |D|\ge
3$,  then $D\equiv 2H-x_1-x_2$ for some points $x_1$
and $x_2$  on $C$.}

By (4.2), the linear system $|2H-(D-x-y)|$ is non-empty for any points $x$ and
$y$ on $C$.
By Riemann-Roch, this implies that $|2H-D|$ is non-empty.

\medskip
{\bf 5. Proof of the theorems}

\ind We only need to prove that, under the hypotheses of theorem 1 and theorem
2, $W_{d-2}(C)$ and
$W_{d-1}(C)$ respectively do not contain any non-zero abelian
varieties. This will follow from the following proposition.

{\pc Proposition} 1. -- {\it Let $C$ be a smooth plane curve of degree $d\ge
4$.  Then:}

{\parindent=1cm\item{\rm (i)}{\it If $d\ge 5$, the variety $W_{d-1}(C)$ does
not contain any abelian
variety of dimension $\ge 2$.
\item {\rm (ii)} If  the variety $W_{d-2}(C)$  contains an
elliptic curve $E$, the inclusion is induced by  a morphism $C\ra E$ of degree
$d-2$ and
$d\le 6$. \item {\rm (iii)} If $d\ge 6$, and if  the variety $W_{d-1}(C)$
contains an
elliptic curve $E$, then $d\le 7$ and the inclusion is induced by  a morphism
$C\ra E$ of
degree $d-1$ or $d-2$.\par}}

{\bf Proof.} If $d=4$, property (ii) follows from [A], theorem
11.2.

\ind We may therefore assume $d\ge 5$. Let $1\le e\le d-1$  and assume that
$W_e(C)$ contains an
abelian variety $A$ of dimension $h>0$. Let $A_2$ be the image of $A\times A$
under the addition
map  $W_e(C)\times W_e(C)\ra W_{2e}(C)$ and let $r$ be the dimension of the
linear system on
$C$ which corresponds to a generic point of  $A_2$. It follows from [A], lemma
8 that $r\ge
h$. We may assume that $A$ is not contained in  $x+W_{e-1}(C)$ for any point
$x$ in $C$. In
this case, the linear system on $C$ that corresponds to a generic point of
$A_2$ is
base-point-free.

\ind Assume first $h\ge 2$. Since $r\ge h$, we get a  family of \bpf linear
systems of
degree $\le 2d-2$ and dimension $\ge 2$ parametrized by $A$, which is an
abelian variety. By
(4.6), this is possible only if $d=5$, $e=d-1=4$ and $h=r=2$. By [A] lemma 14,
the
morphisms $C\ra {\bf P}^2$ which correspond to points of $A_2$ factor through a
fixed
morphism $p:C\ra B$ of degree $n>1$ onto a curve $B$ of genus $\ge h=2$. The
induced
birational morphims $B\ra {\bf P}^2$  then have degree  $8/n$, hence  $n=2$ and
$g(B)=2$.
Let $\sigma$ be the involution associated with the double cover $p$, and let
$H$ be a
hyperplane section of   $C$. Since the embedding of $C$ as a smooth plane
curve is unique (this follows for example from (4.1)), one has
$\sigma^*(H)\equiv H$ hence
$\sigma$ is induced by a projective automorphism $\tau$ of ${\bf P}^2$. By
Riemann-Hurwitz,
$\sigma$ has $6$ fixed points, hence $\tau$ is the symmetry with respect to a
line. But
then, the fixed points of $\sigma$ are the intersection of this line with $C$,
hence
there cannot be $6$ of them since $C$ has degree $5$. Therefore, this case does
not
occur.

\ind This takes care of the case $h\ge 2$ and we now assume $h=1$.

\ind If $r\ge 3$, we get  from (4.6) a non-constant map from the elliptic curve
$A$
into $C^{(2)}$. By [A] theorem 11.2, $C$ is bi-elliptic, hence has an elliptic
curve of
pencils of degree $4$. This contradicts (4.1).

\ind If $r=2$, we get linear systems of dimension $2$ and degree $2e\le 2d-2$.
We get
 $e=d-1$ from (4.2). By (4.6), if moreover  $d\ge 6$,  there are no linear
systems of
degree $2d-2$ and dimension exactly $2$, which is a contradiction.

\ind The only remaining case is $h=r=1$ (except maybe if $d=5$  and $e=4$). The
embedding of the elliptic curve $A$ in $W_e(C)$ is then induced by a morphism
$C\ra A$ of
degree $e$ ([A], lemma 13). The pencils of degree $2$ on $A$ pull back to an
elliptic curve
of \bpf pencils on $C$ of degree $2e$. Fact (4.1) implies $e=d-2$ or $d-1$.

\ind If $e=d-2$, fact (4.2) yields an embedding of $A$
into $W_4(C)$. If $d> 6$, one has $g(C)\ge 8$, and theorem 11 of [A]
implies that $C$ has a morphism of degree $\le 4$ onto an elliptic curve $E$.
But then, the pencils
of degree $2$ on $E$ pull back to an elliptic curve of \bpf pencils on $C$ of
degree $\le 8$, which
contradicts (4.1) (since $8\le 2d-5$.). Therefore, $d\le 6$.

\ind If $e=d-1$ and $d\ge 8$, fact (4.4)  yields an embedding of $A$ into
$W_4(C)$, which we just saw cannot exist. Therefore, $d\le 7$.

\ind This finishes the proof of the proposition.\cqfd
\medskip
{\bf 6. Fermat curves}

\ind Both theorems apply in particular to the Fermat curves $F_d$ defined by
the equation
$X^d+Y^d=Z^d$, at least for $d\ge 8$. For small $d$, the situation is the
following:

\ind $\bullet$ for $d=3$, $4$, $5$ or $7$, it is known  that $J(F_d)(\Q )$,
hence also its
subvarieties $W_e(F_d)(\Q )$ for all $e$, are finite ([F1], [F2]). This of
course implies
both theorems. For $d=4$, Faddeev also shows that in addition to its four
rational points,
$F_4$ has exactly twelve points defined over quadratic fields, and that the
lines through
each of the four rational points of $F_4$ account for {\it all\/} points of
$F_4$ in all
cubic fields.

\ind $\bullet$ for $d=6$, there is a morphism of degree $4$ from
$F_6$  onto the elliptic curve $F_3$. In particular, $W_4(F_6)$ does contain
an elliptic curve and our whole method of proof collapses. However, since
$J(F_3)(\Q )$ is
finite, this does not say anything about the finiteness of $W_4(F_6)(\Q )$. On
the other
hand,   $J(F_6)(\Q )$ is known to be infinite ([F2]). One may use here a
stronger
 theorem of Faltings ([Fa2]), which says that {\it if $X$ is a subvariety of an
abelian variety
$A$, both  defined over a number field $K$, then the set $X(K)$  lies
inside a finite union of $K$--rationally defined translates of abelian
subvarieties of
$A$ contained in $X$}. Consequently,  if  any morphism of
degree $4$ from $F_6$ onto an elliptic curve has image $F_3$, theorem 1 will
hold for
$F_6$. If, in addition,  there are no morphisms of degree $5$ from $F_6$ onto
an elliptic
curve, theorem 2  will hold for $F_6$ .

\ind We mentioned in the introduction that our two theorems remained valid over
any number
field $K$. This applies in particular to Fermat curves for $d\ge 7$ (for
theorem 1) and
$d\ge 8$ (for theorem 2).  Both theorems remain  valid for $F_5$: the
absolutely simple
factors of its Jacobian are surfaces ([KR] theorem 2), and $W_4(F_5)$ cannot
contain an
abelian surface by proposition 1.(i).
However, theorem 1 fails trivially for $F_3$,  for $F_4$ (this
curve has a morphism of degree $2$ onto the elliptic curve $E$ with equation
$U^2W^2+V^4=W^4$, and as soon as   $E(K)$ becomes infinite, so will
$F_4^{(2)}(K)$) and
for $F_6$ (for the same reason, since there is a morphism of degree $4$ from
$F_6$ onto the
elliptic curve $F_3$). As far as $F_7$ is concerned, theorem 1 holds, and
theorem 2 holds if and
only if there are no  morphisms of degree $6$ from $F_7$ onto an elliptic
curve.

\ind We cannot resist  giving a different proof of theorem 1 for Fermat curves
when $d$ is an
odd  {\it prime\/} number $p$ which satisfies $p\equiv 2\!\!\!\!\pmod{3}$ as  a
nice
application of the following result of [DF] (proposition 3.3): {\it for a
complex projective
curve $C$ of genus $g$, the variety $W_d(C)$ cannot contain an abelian variety
of dimension
$> d/2$ for $d<g$.} In fact, it is known in this case that the absolutely
simple factors of
$J(F_p)$ are all of dimension ${p-1\over 2}$ ([KR], theorem 2), so that
$W_{p-2}(F_p)$
cannot contain any non-zero abelian variety.

\ind It would be very interesting to know specifically which points constitute
the finite
sets of algebraic points in the theorems for the Fermat curves.  Of course,
this
extends the already difficult question, posed by Fermat, of showing that there
are only
three rational points if $d$ is odd and four if $d$ is even.

\ind Assume again that $d$ is an odd prime number $p$. The easiest way to
produce  algebraic
points of degree $\le p-3$ is to take the other points of intersection of
$F_p$ with the line through the three known rational points. In the affine
patch where $Z=1$, this is just the line $y=1-x$, and  $F_p$ is defined by
$x^p+y^p=1$. The $x$-coordinates of these other points of intersection are
then just the roots of $${x^p+(1-x)^p-1 \over x(x-1)} = 0\ .$$
One sees,
by considering the equation $x^p+(1-x)^p-1 = 0$ and its derivative,
 that  the factor $x^2-x+1$ always occurs with multiplicity
one or two depending on whether $p$ is $5$ or $1\!\!\!\!\pmod{6}$ respectively,
so we
obtain $x$-coordinates $\eta$, and $\eta^{-1}$, where $\eta$ is a primitive
sixth root of unity. The other factor, of degree $p-5$ or $p-7$, is irreducible
over $\Q$ for $p \le 101$ (checked using MAPLE), but the authors do not know
if this is always the case. Also, the authors do not know of any other
points on $F_p$ of degree $\le p-2$, \ie which do not lie on the line
 $y = 1 - x$. Gross and Rohrlich show that this line accounts for all the
 points of degree $\le (p-1)/2$ on $F_p$ for the
primes $p = 3$, $5$, $7$ and $11$ (see [GR].)

\ind It is interesting to note that the linear equivalence class described
above
produces the only known points of infinite order on the Mordell-Weil group
$J(F_p)(\Q )$.  More specifically, Gross and Rohrlich take the conjugate
quadratic points $P = (\eta,\eta^{-1},1)$, and $\overline{P} = (\eta^{-1},
\eta,1)$,
and form the divisor $P + \overline{P} -2\infty$ on $\Div ^0(F_p)$. Then, they
show that for $p>7$, the linear equivalence class in
$J(F_p)$ of the divisor  $P + \overline{P}
-2\infty$ represents a point of infinite order.

\ind Finally, one would like to have, if not a complete description, at least
an upper bound on the cardinality of the finite sets of algebraic points
in  theorem 1. It is natural to begin by trying to bound the number
of $\Q$--rational points on $F_p$. The greatest success in this regard is
Kummer's proof of Fermat's Last Theorem for $p$ a regular prime (see [W].)
For general $p$, all bounds depend on the rank of $J(F_p)(\Q )$.  One approach
is
to come up with bounds which are
exponential in the Mordell-Weil rank as in Bombieri's version of Faltings'
Theorem ([B].)  Another approach is to use Coleman's ``Effective
Chabauty,'' applying his theory of $p$-adic abelian integrals ([C1], [C2].) In
this
case, one needs to know that the rank of the Mordell-Weil group  $J(F_p)(\Q )$
is
less than its dimension (\ie the genus of the curve.)  This
is  known to hold in the case when $p$ is regular ([F3].)
However, McCallum has shown that this would hold for all $p$ if one has a
certain bound on the ideal class group $\Cl$ of $\Q (\zeta)$,
where $\zeta$ is a primitive $p^{th}$ root of unity.  In particular, if $\Cl
[p]$
denotes the subgroup of $\Cl$ of elements killed by $p$, then he shows that
$$ \rank_{\bf Z}J_s(\Q ) \le {p-7\over 4} + 2\rank_{{\bf Z} /p{\bf Z} }(\Cl
[p]).
$$ He then goes on to show that if $\rank_{\Z /p\Z }(\Cl [p]) < {p+5\over 8}$,
then
the number of $\Q$--rational points on $F_p$ is $\le 2p-3$ ([Mc].)
The second author has begun to apply Coleman's theory to the symmetric products
of
curves in his Ph.D. thesis (to appear.)

 \saut
\saut  \centerline {\pc References}

\saut

\hangindent=1cm
[A]	Abramovich, D., {\it Subvarieties of abelian varieties and of Jacobians of
curves\/}. Ph.D. Thesis, Harvard University, 1991.

\hangindent=1cm
[ACGH]  Arbarello, E., Cornalba, M., Griffiths, P., Harris, J.,
       {\it Geometry of algebraic curves. I. Grundlehren 267\/},
Springer-Verlag, New York, 1985.

\hangindent=1cm
[AH]	Abramovich, D.,
Harris, J., {\it Abelian varieties and curves on  $W_d(C)$\/}, Comp. Math. 78
(1991), 227--238.

\hangindent=1cm
[B]	 Bombieri, E., {\it The Mordell Conjecture Revisited\/},
         Preprint, 1991.

\hangindent=1cm
[C1] Coleman, R., {\it Torsion points on curves and $p$-adic
        abelian integrals\/}, Ann. of Math. 121 (1985),
        111--168.

\hangindent=1cm
[C2] Coleman, R., {\it Effective Chabauty
    \/}, Duke Math. J. 52 (1985), 765--770.

\hangindent=1cm
[CK] Coppens, M., Kato, T., {\it The gonality of smooth curves with plane
models\/},
Manuscripta Math.  70    (1990),  5--25, and {\it Correction to: ``the gonality
of smooth
curves with plane models"\/}, Manuscripta Math.  71    (1991),  337--338.

\hangindent=1cm
[DF] Debarre, O., Fahlaoui, R., {\it Abelian Varieties In  $W_d^r(C)$  And
Points
Of Bounded Degrees On Algebraic Curves\/}. To appear in Comp. Math.

\hangindent=1cm
[F1] Faddeev, D.K., {\it Group of divisor classes on the curve defined by the
equation
$x^4+y^4=1$\/}, Dokl.
Akad. Nauk. SSSR 134 (1960), 776--778 and Sov. Math. Dokl. 1 (1960),
1149--1151.

\hangindent=1cm
[F2] Faddeev, D.K., {\it On the divisor class group of some algebraic
curves\/}, Dokl.
Akad. Nauk. SSSR 136 (1960), 296--298 and Sov. Math. Dokl. 2 (1961), 67--69.

\hangindent=1cm
[F3] Faddeev, D.K., {\it Invariants of divisor classes for the
        curves $x^k(1-x) = y^l$ in an $l$-adic cyclotomic field\/} (in
Russian), Trudy Math.
        Inst. Steklov  64    (1961),  284--293.

\hangindent=1cm
[Fa1] Faltings, G., {\it Diophantine approximation on abelian
        varieties\/}, Ann. of Math.  133   (1991),
549--576.

\hangindent=1cm
[Fa2] Faltings, G., {\it The general case of S. Lang's conjecture\/}, to
appear.

\hangindent=1cm
[GR] Gross, B.H., Rohrlich, D.E., {\it Some results on the
        Mordell-Weil group of the Jacobian of a quotient of the Fermat
curve\/}, Invent. Math.  93 (1978),  637--666.

\hangindent=1cm
[KR] Koblitz, N., Rohrlich, D.E., {\it Simple Factors in the
        Jacobian of a Fermat Curve}, Canadian J. Math.  30
         (1978),   1183--1205.

%\hangindent=1cm
%[M] Mumford, D., {\it Prym Varieties, I\/}, in ``Contributions to Analysis",
%Academic Press, New York, 1974.

\hangindent=1cm
[Mc] McCallum, W.G., {\it The Arithmetic of Fermat Curves\/}, to appear in
Math.
Ann.

[W] Washington, L.C., {\it Introduction to Cyclotomic Fields\/}, Springer
Verlag,
New York, 1982.
\hangindent=1cm

\bye